\documentclass[11pt,a4paper,reqno]{amsart}
\usepackage{mathrsfs}
\usepackage{amsfonts}
\usepackage{amsthm}
\usepackage{lmodern}
\usepackage[textwidth=16cm, hmarginratio=1:1]{geometry}
\usepackage{amsmath,amssymb,mathtools}
\usepackage{microtype,xcolor}
\usepackage[
colorlinks=true,
linkcolor=blue!45!black,
citecolor=blue!45!black,
urlcolor=blue!45!black
]{hyperref}

\hypersetup{
pdftitle={Normality of marginal inverse subsemigroups in inverse semigroups}}

\newtheorem{theorem}{Theorem}[section]
\newtheorem{corollary}[theorem]{Corollary}
\newtheorem{lemma}[theorem]{Lemma}
\newtheorem{proposition}[theorem]{Proposition}

\newtheorem{example}[theorem]{Example}

\newtheorem{problem}[theorem]{Problem}

\newcommand{\C}{\mathbb C}
\newcommand{\Sn}{S_{m,n}}
\newcommand{\Mn}{M_{x^n}}
\newcommand{\E}{E}
\newcommand{\cent}[1]{C_{#1}\bigl(E(#1)\bigr)}
\newcommand{\id}{\operatorname{id}}
\DeclareMathOperator{\rank}{rank}
\DeclareMathOperator{\Fix}{Fix}
\DeclareMathOperator{\dom}{dom}
\DeclareMathOperator{\ran}{ran}
\numberwithin{equation}{section}

\title[Normality of marginal inverse subsemigroups in inverse semigroups]{Normality of marginal inverse subsemigroups in inverse semigroups}
\author{Bing Duan}
\address{Bing Duan: School of Mathematics and Statistics, Lanzhou University, Lanzhou 730000, P. R. China.}
\email{duanbing@lzu.edu.cn}
\author{Jizhi Zhang}
\address{Jizhi Zhang: School of Mathematics and Statistics, Lanzhou University, Lanzhou 730000, P. R. China.}
\email{zhangjzh2024@lzu.edu.cn}
\date{}

\begin{document}

\maketitle

\begin{abstract}
For each pair of integers $m\ge1$ and $n\ge2$, we construct a finite inverse semigroup $\Sn \subseteq M_m(\mathbb C)$ and explicitly determine its marginal inverse subsemigroup $\Mn(\Sn)$ associated with the word $x^n$. We prove that the three conditions $\sim_{x^n}=\equiv_{x^n}$ on $S_{m,n}$, $\Mn(\Sn)\trianglelefteq\Sn$, and $m\le n$ are equivalent. As a consequence, we answer Problem~11 of Ara\'ujo, Ara\'ujo and Kinyon in \cite{AAK}. Moreover, for a word $\omega$ in any inverse semigroup $S$, we prove that $\sim_\omega=\equiv_\omega$  if and only if $M_\omega(S)\trianglelefteq S$ and $\sim_\omega$ is invariant under inversion.

\hspace{0.15cm}


\hspace{0.15cm}

\end{abstract}

\section{Introduction}

A semigroup $S$ is called an inverse semigroup if for each $x\in S$, there exists a unique element $x^{-1}\in S$ such that
\[
xx^{-1}x=x \qquad\text{and}\qquad x^{-1}xx^{-1}=x^{-1}.
\]
Throughout, let $S$ be an inverse semigroup, and denote its set of idempotents by $E(S)$. An inverse subsemigroup $N$ of $S$ is full if $E(S)\subseteq N$, and normal if it is full and $sNs^{-1}\subseteq N$ for each $s\in S$. The latter condition is also called self-conjugacy. We write $N\trianglelefteq S$ for normality. For standard background on inverse semigroups, see \cite{Howie,Lawson}.

Marginal subgroups, introduced by P.~Hall \cite{Hall}, are characteristic, and hence normal, subgroups of a group. Ara\'ujo, Ara\'ujo and Kinyon \cite{AAK} extended this notion to inverse semigroups. Let $\omega(x_1,\ldots,x_d)$ be a word in $S$. An element $a\in S$ is called marginal for $\omega$ if, for each $i\in\{1,\ldots,d\}$, there exist idempotents $e_i,f_i\in E(S)$ such that, for every $\boldsymbol{x}=(x_1,\ldots,x_d)\in S^d$,
\begin{align}
\omega(x_1,\ldots,ax_i,\ldots,x_d) & = \omega(x_1,\ldots,e_ix_i,\ldots,x_d), \label{eq:def-left}\\
\omega(x_1,\ldots,x_ia,\ldots,x_d) & = \omega(x_1,\ldots,x_if_i,\ldots,x_d). \label{eq:def-right}
\end{align}
The set of all marginal elements for $\omega$ is denoted by $M_\omega(S)$. By \cite[Theorem~1]{AAK}, $M_\omega(S)$ is a full inverse subsemigroup of $S$. Moreover, \cite[Lemma~4]{AAK} shows that every $a\in M_\omega(S)$ is completely regular and that the witnessing idempotents are uniquely determined by $e_i=f_i=a^0:=aa^{-1}=a^{-1}a$ for all $i$. By \cite[Lemma~6]{AAK}, every marginal element centralizes $E(S)$. Consequently,
\begin{equation}\label{eq:bounds}
E(S)\subseteq M_\omega(S)\subseteq \cent S,
\end{equation}
where $\cent S$ denotes the centralizer of $E(S)$ in $S$.

Following \cite[Section 5]{AAK}, we define the marginal equivalence relation $\sim_\omega$ on $S$ as follows: for $a,b\in S$,
we write $a\sim_\omega b$ if, for every $i\in\{1,\ldots,d\}$ and every
$\boldsymbol{x}=(x_1,\ldots,x_d)\in S^d$,
\[
\begin{aligned}
\omega(x_1,\ldots,ax_i,\ldots,x_d) &=\omega(x_1,\ldots,bx_i,\ldots,x_d),\\
\omega(x_1,\ldots,x_i a,\ldots,x_d) &=\omega(x_1,\ldots,x_i b,\ldots,x_d).
\end{aligned}
\]
The marginal congruence $\equiv_\omega$ is the inverse semigroup congruence generated by $\sim_\omega$. For the power word $\omega(x)=x^n$, where $n\ge1$, the marginal
equivalence relation takes the form
\[
a\sim_{x^n}b \quad\Longleftrightarrow\quad (ax)^n=(bx)^n \ \text{and}\ (xa)^n=(xb)^n \quad\text{for all }x\in S.
\]
The marginal inverse subsemigroup associated with $x^n$ is given by
\begin{equation}\label{eq:power-margin}
\Mn(S)=
\left\{
a\in S \;\middle|\;
\begin{aligned}
(ax)^n&=(aa^{-1}x)^n,\\
(xa)^n&=(xaa^{-1})^n
\end{aligned}
\quad\text{for all }x\in S
\right\}.
\end{equation}

Recall that the metacenter of an inverse semigroup $S$ is
\[
Z(S)=\{a\in S: aa^{-1}xa=axa^{-1}a \text{ for all } x\in S\}.
\]
The metacenter is known to be a normal inverse subsemigroup of $S$. For the commutator word, Ara\'ujo, Ara\'ujo, and Kinyon proved in \cite[Theorem~2]{AAK} that $M_\omega(S)$ coincides with the metacenter of $S$ and is therefore normal.  For a general inverse semigroup word $\omega$, they posed the following question:

\begin{problem}[\cite{AAK}, Problem~11] \label{AAK-Problem11}
Let $\omega(x_1,\ldots,x_n)$ be a word in $S$. Do the marginal congruence $\equiv_\omega$ and the marginal equivalence relation $\sim_\omega$ coincide? Equivalently, is the marginal inverse subsemigroup $M_\omega(S)$ normal in $S$?
\end{problem}

For each pair of integers $m\geq 1$ and $n\geq 2$, we construct a finite inverse semigroup $\Sn\subseteq M_m(\mathbb C)$ and explicitly determine its marginal inverse subsemigroup $\Mn(\Sn)$ for the power word $x^n$, see Theorem \ref{thm:margin}. We prove in Theorem \ref{thm:normality} that
\[
\sim_{x^n}=\equiv_{x^n}\ \text{on }\Sn \quad\Longleftrightarrow\quad \Mn(\Sn)\trianglelefteq\Sn \quad\Longleftrightarrow\quad m\leq n.
\]
Thus, for $m>n\geq 2$, we have $\sim_{x^n}\subsetneq\equiv_{x^n}$ on $\Sn$, and the margin $\Mn(\Sn)$ is not normal in $\Sn$.

We also clarify the relation between the marginal equivalence relation $\sim_\omega$ and the marginal congruence $\equiv_\omega$ in Problem~\ref{AAK-Problem11}. For an arbitrary inverse semigroup word $\omega$ in $S$, we prove in Theorem \ref{thm:general-equivalence} that $\sim_\omega$ is a congruence if and only if $M_\omega(S)\trianglelefteq S$ and $\sim_\omega$ is invariant under inversion. Invariance under inversion holds automatically for power words, so the two formulations are equivalent in that case, see Lemma \ref{lem:symmetric-inversion} and Corollary \ref{cor:symmetric-normality}. We give an example for a general inverse semigroup word in which the margin is normal, while marginal equivalence fails to be invariant under inversion, see Propositions \ref{prop:omega-margin} and \ref{prop:omega-inversion}. This example shows that the two formulations of Problem~\ref{AAK-Problem11} are not equivalent in general, since $M_\omega(S)$ can be normal even when $\sim_\omega\subsetneq\equiv_\omega$.

The organization of this paper is as follows. In Section~\ref{sec:structures}, we define an inverse semigroup $\Sn  \subseteq  M_m(\mathbb C)$ and explicitly determine its marginal inverse subsemigroup $\Mn(\Sn)$. In Section~\ref{sec:normality}, we prove Theorem \ref{thm:normality} and characterize the marginal elements that have a conjugate outside $\Mn(\Sn)$ when $m>n$. In Section~\ref{sec:equivalence}, we study the relation between normality and marginal congruences, and prove that, for an inverse semigroup word $\omega$,  $\sim_\omega=\equiv_\omega$ if and only if $M_\omega(S)\trianglelefteq S$ and  $\sim_\omega$ is invariant under inversion. In Section~\ref{sec:generality}, we give two situations in which the bounds \eqref{eq:bounds} become equalities.

\section{Marginal inverse subsemigroup $\Mn(\Sn)$ of $\Sn$} \label{sec:structures}

Fix integers $m\geq 1$ and $n\geq 2$, and let
\[
[m]=\{1,\ldots,m\},\quad \mu_n=\{\zeta\in\C:\zeta^n=1\},\quad L=\min\{m,n\}.
\]
For nonempty subsets $I,J\subseteq[m]$ of the same cardinality $r$, we list their elements in increasing order as $I=\{i_1<\cdots<i_r\}$ and $J=\{j_1<\cdots<j_r\}$. For $\zeta\in\mu_n$, let
\[
A(I,J;\zeta) =\zeta E_{i_1j_1}+\sum_{k=2}^{r} E_{i_kj_k},
\]
where $E_{ij}$ denotes the standard matrix units in $M_m(\C)$. We then set
\begin{equation*}
\Sn=\{0\}\cup \bigl\{A(I,J;\zeta)\bigm| \, 1\leq |I|=|J|\leq L,\ \zeta\in\mu_n \bigr\}.
\end{equation*}
For simplicity of notation, for each nonempty subset $I\subseteq[m]$ and each $\zeta\in\mu_n$, let
\[
d_I(\zeta)=A(I,I;\zeta),\quad P_I=d_I(1)=\sum_{i\in I}E_{ii}.
\]
We call $I$ the support of the diagonal matrix $d_I(\zeta)$.

For each nonzero $A=A(I,J;\zeta)\in\Sn$, we define its support map $\phi_A$ to be the increasing partial bijection of $[m]$ with domain $J$ and range $I$, given by
\[
\phi_A(j_k)=i_k \qquad (1\leq k\leq r).
\]
We write $D_A=\dom\phi_A=J$ and $R_A=\ran\phi_A=I$. Then
$$
\operatorname{rank}(\phi_A)=\operatorname{rank}(A)=r\leq L.
$$
For the zero matrix, we take $\phi_0$ to be the empty map.


For $A,B\in\Sn$, each nonzero entry of $AB$ comes from a unique matching pair of nonzero entries of $A$ and $B$. Thus the support map of the product satisfies $\phi_{AB}=\phi_A\circ\phi_B$. Once we establish that $\Sn$ is an inverse semigroup, this identity shows that the map
$$
\pi\colon \Sn\longrightarrow\mathcal I_m,\qquad A\longmapsto\phi_A,
$$
is an inverse semigroup homomorphism, where $\mathcal I_m$ is the symmetric inverse semigroup of all partial bijections of $[m]$. Moreover, for each $A\in\Sn$,
\begin{equation}\label{eq:diagonal}
A\text{ is diagonal}
\quad\Longleftrightarrow\quad
\phi_A=\id_{D_A}
\quad\Longleftrightarrow\quad
D_A=R_A.
\end{equation}
Indeed, if $D_A=R_A$, then $\phi_A$ is an increasing bijection from the finite ordered set $D_A$ to itself and hence is the identity.

\begin{proposition}\label{prop:inverse}
The set $\Sn$ is a finite inverse semigroup under matrix multiplication. The inverse of $A\in\Sn$ is its conjugate transpose $A^*=\overline{A}^{\,T}$. Moreover,
\begin{align}
|\Sn| &=1+n\sum_{r=1}^{L}\binom{m}{r}^{2},\label{eq:cardinality}\\
\E(\Sn) &=\{0\}\cup\{P_I:I\subseteq[m],\ 1\leq |I|\leq L\}. \label{eq:idempotents}
\end{align}
\end{proposition}

\begin{proof}
Let $A=A(I,J;\zeta)$ and $B=A(K,H;\eta)$ be nonzero elements of $\Sn$. If $J\cap K=\varnothing$, then $AB=0$. Otherwise, each $t\in J\cap K$ contributes a nonzero entry in row $\phi_A(t)$ and column
$\phi_B^{-1}(t)$. Both indices increase with $t$, so $AB$ has an increasing support map of rank $|J\cap K|\leq L$.

The entry corresponding to $t$ can receive a factor of $\zeta$ only when $t=\min J$ and a factor of $\eta$ only when $t=\min K$. Either condition forces $t=\min(J\cap K)$. Thus the first nonzero entry of
$AB$ has coefficient in $\{1,\zeta,\eta,\zeta\eta\}\subseteq\mu_n$, and all subsequent nonzero entries have coefficient $1$. Hence $AB\in\Sn$. Products involving the zero matrix also lie in $\Sn$.

Since $\overline{\zeta}=\zeta^{-1}$, we have
$$
A^*=A(J,I;\zeta^{-1})\in\Sn,\qquad
AA^*=P_I,\qquad A^*A=P_J.
$$
It follows that $AA^*A=A$ and $A^*AA^*=A^*$, with both identities also holding for $A=0$. Thus $\Sn$ is regular.

If $A$ is a nonzero idempotent, then $\phi_A^2=\phi_A$, so $\phi_A$ is a partial identity. By \eqref{eq:diagonal}, $A=d_I(\zeta)$ for some $I$ and $\zeta$. The identity $A^2=A$ gives $\zeta^2=\zeta$, hence $\zeta=1$ and $A=P_I$. Conversely, $0$ and all the matrices $P_I$ in \eqref{eq:idempotents} are idempotent. This proves \eqref{eq:idempotents}. These idempotents commute, so regularity implies that $\Sn$ is an inverse semigroup.

For each $1\leq r\leq L$, there are $\binom{m}{r}^2$ choices of $(I,J)$ with $|I|=|J|=r$ and $n$ choices of $\zeta$. The nonzero rows and columns determine $I$ and $J$, and the first nonzero entry determines $\zeta$. Counting these distinct matrices together with $0$ gives \eqref{eq:cardinality}.
\end{proof}

\begin{lemma}\label{lem:paths}
Let $\phi$ be an increasing partial bijection of $[m]$ with $\rank\phi\leq n$. If at least one of the conditions
\[
\rank\phi<n,\quad  \{1,m\}\subseteq\dom\phi,\quad \{1,m\}\subseteq\ran\phi
\]
holds, then $\phi^n=\id_{\Fix(\phi)}$.
\end{lemma}

\begin{proof}
Every fixed point of $\phi$ lies in the domain of $\phi^n$ and is fixed by $\phi^n$. Suppose that $\phi^n(t)$ is defined for some $t\notin\Fix(\phi)$, and put $t_k=\phi^k(t)$ for $0\leq k\leq n$, with $t_0=t$. Since $\phi$ is increasing and $t_1\neq t_0$, the sequence $t_0,\ldots,t_n$ is strictly monotone. The rank bound then gives
\[
\rank\phi=n,\quad  \dom\phi=\{t_0,\ldots,t_{n-1}\},\quad  \ran\phi=\{t_1,\ldots,t_n\}.
\]
If the sequence is increasing, every domain element is less than $t_n\leq m$ and every range element is greater than $t_0\geq 1$. Thus $m\notin\dom\phi$ and $1\notin\ran\phi$. If the sequence is
decreasing, every domain element is greater than $t_n\geq 1$ and every range element is less than $t_0\leq m$. Thus $1\notin\dom\phi$ and $m\notin\ran\phi$.

In either case, none of the three conditions holds, a contradiction. Hence $\dom(\phi^n)=\Fix(\phi)$, and $\phi^n=\id_{\Fix(\phi)}$.
\end{proof}

\begin{corollary}\label{cor:matrix-power}
Let $Y\in\Sn$. If its support map $\phi_Y$ satisfies any of the three conditions in Lemma~\ref{lem:paths}, then $Y^n=P_{\Fix(\phi_Y)},$  where $P_\varnothing=0$.
\end{corollary}

\begin{proof}
By multiplicativity of the support map and Lemma~\ref{lem:paths},
$$
\phi_{Y^n}=\phi_Y^n=\id_{\Fix(\phi_Y)}.
$$
Thus the nonzero entries of $Y^n$ occur precisely at the diagonal
positions indexed by $\Fix(\phi_Y)$. For each $i\in\Fix(\phi_Y)$,
the only nonzero entry in column $i$ of $Y$ is $Y_{ii}\in\mu_n$.
Consequently, $(Y^n)_{ii}=(Y_{ii})^n=1$. This gives
$Y^n=P_{\Fix(\phi_Y)}$.
\end{proof}

\begin{theorem}\label{thm:margin}
For $m\geq 1$ and $n\geq 2$,
\begin{equation}\label{eq:margin-formula}
\Mn(\Sn)=\{0\}\cup
\left\{
d_I(\zeta)\;\middle|\;
\begin{array}{l}
I\subseteq[m],\quad 1\leq |I|\leq L,\quad \zeta\in\mu_n,\\
\zeta=1\ \text{or}\ |I|<n\ \text{or}\ \{1,m\}\subseteq I
\end{array}
\right\}.
\end{equation}
\end{theorem}

\begin{proof}
Let $A\in\Mn(\Sn)$ be nonzero. Taking $X=AA^*$ in the right-hand
equality of \eqref{eq:power-margin} gives
$A^n=AA^*=P_{R_A}$, since $AA^*A=A$ and $AA^*$ is idempotent.
Thus $\phi_A^n=\id_{R_A}$. Since $\dom(\phi_A^n)\subseteq D_A$, we
obtain $R_A\subseteq D_A$. As $R_A$ and $D_A$ have the same
cardinality, they are equal. By \eqref{eq:diagonal}, $A$ is diagonal.

The zero matrix is marginal. For a nonzero diagonal element, write
$A=d_I(\zeta)$, where $I=\{i_1<\cdots<i_r\}$ and $AA^*=P_I$.
If $\zeta=1$, then $A=P_I$ is idempotent and hence marginal.

Suppose that $r<n$ or $\{1,m\}\subseteq I$, and let $X\in\Sn$.
Since $\phi_A=\phi_{P_I}=\id_I$, the matrices $AX$ and $P_IX$ have
the same support map, with range $I\cap R_X$. Likewise, $XA$ and
$XP_I$ have the same support map, with domain $I\cap D_X$.
Each of these support maps has rank at most $r\leq n$.
If the first has rank $n$, then $I\cap R_X=I$ and $I$ contains both
endpoints. If the second has rank $n$, then $I\cap D_X=I$ and again
$I$ contains both endpoints. Thus Corollary~\ref{cor:matrix-power}
applies to all four products, giving
$$
(AX)^n=(P_IX)^n=P_{\Fix(\phi_{AX})},
\qquad
(XA)^n=(XP_I)^n=P_{\Fix(\phi_{XA})}.
$$
By \eqref{eq:power-margin}, $A$ is marginal.

It remains to consider the case $\zeta\neq 1$, $r=n$, and
$\{1,m\}\nsubseteq I$. If $i_1>1$, put $i_0=i_1-1$ and take
$$
Y=\sum_{k=1}^{n}E_{i_k i_{k-1}}\in\Sn.
$$
Both $AY$ and $P_IY$ have support map given by the chain
$i_0\mapsto i_1\mapsto\cdots\mapsto i_n$. The first arrow has
coefficient $\zeta$ in $AY$ and coefficient $1$ in $P_IY$, while
all other arrows have coefficient $1$ in both products. Hence
$$
(AY)^n=\zeta E_{i_n i_0}\neq E_{i_n i_0}=(P_IY)^n,
$$
so the left-hand equality of \eqref{eq:power-margin} fails.

If $i_n<m$, put $i_{n+1}=i_n+1$ and take
$$
Y=\sum_{k=1}^{n}E_{i_{k+1}i_k}\in\Sn.
$$
The support maps of $YA$ and $YP_I$ are given by the chain
$i_1\mapsto i_2\mapsto\cdots\mapsto i_{n+1}$. Again, only the
first arrow has coefficient $\zeta$ in $YA$, while every arrow
has coefficient $1$ in $YP_I$. Thus
$$
(YA)^n=\zeta E_{i_{n+1}i_1}\neq E_{i_{n+1}i_1}=(YP_I)^n,
$$
so the right-hand equality of \eqref{eq:power-margin} fails.
Since $I$ omits at least one endpoint, at least one of these cases
occurs. Therefore $A$ is not marginal.
\end{proof}

\begin{corollary}\label{cor:count}
For $m\geq 1$ and $n\geq 2$,
$$
|\Mn(\Sn)|
=1+\sum_{r=1}^{L}\binom{m}{r}
 +(n-1)\left(
     \sum_{r=1}^{\min\{m,n-1\}}\binom{m}{r}+h_{m,n}
   \right),
$$
where
$$
h_{m,n}=
\begin{cases}
\displaystyle\binom{m-2}{n-2}, & m\geq n,\\
0, & m<n.
\end{cases}
$$
In particular, for $m=n+1$,
$$
|S_{n+1,n}|=n\binom{2n+2}{n+1}-2n+1,
\qquad
|\Mn(S_{n+1,n})|=n2^{n+1}-4n+3.
$$
\end{corollary}

\begin{proof}
The zero matrix and the matrices $P_I$ with $1\leq |I|\leq L$
contribute $1+\sum_{r=1}^{L}\binom{m}{r}$ idempotents. By
Theorem~\ref{thm:margin}, each nonempty support of size less than $n$
contributes $n-1$ further elements, one for each
$\zeta\in\mu_n\setminus\{1\}$. A support of size $n$ contributes
these elements exactly when it contains both $1$ and $m$.
For $m\geq n$, the remaining $n-2$ points can be chosen from the
$m-2$ interior points in $\binom{m-2}{n-2}$ ways. For $m<n$,
there are no supports of size $n$.

Now let $m=n+1$, so $L=m-1$ and $h_{m,n}=m-2$. By Vandermonde's
identity and the binomial theorem,
$$
\sum_{r=1}^{m-1}\binom{m}{r}^{2}=\binom{2m}{m}-2,\qquad
1+\sum_{r=1}^{m-1}\binom{m}{r}=2^m-1,\qquad
\sum_{r=1}^{m-2}\binom{m}{r}+h_{m,n}=2^m-4.
$$
Substituting these identities into \eqref{eq:cardinality} and the
formula for $|\Mn(\Sn)|$ gives the two asserted values.
\end{proof}

\section{Normality and marginal congruences}\label{sec:normality}

\begin{theorem}\label{thm:normality}
Let $m\ge1$ and $n\ge2$. Then
\[
\begin{aligned}
&\sim_{x^n}=\equiv_{x^n} \text{ on }\Sn \quad\Longleftrightarrow\quad \Mn(\Sn)\trianglelefteq\Sn \quad\Longleftrightarrow\quad m\le n.
\end{aligned}
\]
\end{theorem}

\begin{proof}
Suppose first that $m\le n$. If a support has rank $n$, then $m=n$ and its domain equals $[m]$. Theorem~\ref{thm:margin} therefore says that $\Mn(\Sn)$ consists of all diagonal elements of $\Sn$. Conjugation by an element in $\Sn$ sends a diagonal matrix to a diagonal matrix: at the level of support maps,
\[
\phi_{XAX^*}=\phi_X\circ\id_{D_A}\circ\phi_X^{-1}
=\id_{\phi_X(D_A\cap D_X)}.
\]
Thus $\Mn(\Sn)$ is self-conjugate and is normal.

We now prove that $A\sim_{x^n}B$ if and only if $\phi_A=\phi_B$. Since $m\le n$, every $A\in\Sn$ satisfies either $\operatorname{rank}(\phi_A)<n$, or $m=n$ and $D_A=R_A=[m]$. Thus the conditions of Lemma~\ref{lem:paths} are satisfied, giving $A^n=P_{\operatorname{Fix}(\phi_A)}$. Suppose first that $\phi_A=\phi_B$. For every $X\in\Sn$, we have
\[
\phi_{AX}=\phi_A\circ\phi_X=\phi_B\circ\phi_X=\phi_{BX},
\]
and hence $(AX)^n=(BX)^n$. Similarly, $\phi_{XA}=\phi_{XB}$ implies $(XA)^n=(XB)^n$. Therefore $A\sim_{x^n}B$. Conversely, suppose that $\phi_A\ne\phi_B$. Interchanging $A$ and $B$ if necessary, we may choose $j$ such that $\phi_A(j)=i$, while $\phi_B(j)\ne i$ or $\phi_B(j)$ is undefined. Then $AE_{ji}=\xi E_{ii}$ for some $\xi\in\mu_n$, so $(AE_{ji})^n=E_{ii}$. On the other hand, $BE_{ji}$ is either zero or a scalar multiple of a non-diagonal matrix unit, and hence $(BE_{ji})^n=0$ for $n\ge2$. Consequently, $A\not\sim_{x^n}B$. We have therefore established that $A\sim_{x^n}B$ if and only if $\phi_A=\phi_B$. Since $\pi:\Sn\to\mathcal I_m,\ A\mapsto\phi_A$, is an inverse semigroup homomorphism, equality of support maps defines a congruence on $\Sn$. Thus $\sim_{x^n}$ is a congruence, and consequently $\sim_{x^n}=\equiv_{x^n}$.

Now suppose that $m>n$. Choose $I=\{i_1<\cdots<i_n\},\,i_1=1,\,i_n=m,\,\zeta\in\mu_n\setminus\{1\}$, and let $A=d_I(\zeta)$. By Theorem~\ref{thm:margin}, $A\in\Mn(\Sn)$. Let $\ell$ be the largest integer such that $i_k=k$ for $1\le k\le\ell$. Since $i_n=m>n$, we have $\ell<n$ and $i_{\ell+1}\ge\ell+2$. The matrix
\begin{equation}\label{eq:conjugator}
X=\sum_{k=1}^{\ell}E_{k+1,k}+\sum_{k=\ell+1}^{n}E_{i_ki_k}
\end{equation}
belongs to $\Sn$. Direct multiplication gives
\begin{equation*}
B:=XAX^*=\zeta E_{22}+\sum_{k=2}^{\ell}E_{k+1,k+1}
+\sum_{k=\ell+1}^{n}E_{i_ki_k}.
\end{equation*}
This diagonal matrix has rank $n$, its least support element is $2$, and its first coefficient is not $1$. Theorem~\ref{thm:margin} gives $B\notin\Mn(\Sn)$, which proves that the margin is not normal.

For the congruence assertion, marginality gives $A\sim_{x^n}AA^*$, and hence
\[
B=XAX^*\equiv_{x^n}XAA^*X^*=BB^*.
\]
Since $B\notin\Mn(\Sn)$, \eqref{eq:power-margin} gives $B\not\sim_{x^n}BB^*$. Thus $\sim_{x^n}\subsetneq\equiv_{x^n}$.
\end{proof}

\begin{corollary}\label{cor:monoids}
For every $n\ge2$, there is a finite inverse monoid whose margin for $x^n$ is not normal and whose marginal equivalence is not a congruence.
\end{corollary}

\begin{proof}
Take $m>n$ and adjoin the identity matrix $I_m$ to $\Sn$. Every element of $\Mn(\Sn)$ is diagonal with coefficients in $\mu_n$, so $A^n=AA^*$. Consequently the additional test $X=I_m$ in \eqref{eq:power-margin} imposes no new restriction on old marginal elements. Old non-marginal elements retain their separating witnesses. Therefore
\[
\Mn(\Sn\cup\{I_m\})=\Mn(\Sn)\cup\{I_m\}.
\]
The conjugate and the separating pair in Theorem~\ref{thm:normality} remain valid.
\end{proof}

\begin{example}\label{ex:37}
Let $m=3,n=2$. Then
\[
|S_{3,2}|=37,\qquad |M_{x^2}(S_{3,2})|=11.
\]
Explicitly, the margin consists of $0$, the six nonzero idempotents of ranks at most $2$, the three matrices $-E_{ii}$, and $A=-E_{11}+E_{33}$. With $X=E_{21}+E_{33}$, one obtains $B=XAX^*=-E_{22}+E_{33}\notin M_{x^2}(S_{3,2})$: for $Y=E_{21}+E_{32}$,
\[
(BY)^2=-E_{31},\qquad ((BB^*)Y)^2=E_{31}.
\]
Adjoining $I_3$ gives an inverse monoid of order $38$.
\end{example}

The precise location of the margin within the idempotent centralizer also explains the failure of normality.

\begin{proposition}\label{prop:layers}
The centralizer $\cent{\Sn}$ consists of all diagonal elements. For $j=1,2,3$, let $T_j$ consist of the non-idempotent diagonal matrices $d_I(\zeta)$, $\zeta\ne1$, satisfying respectively
\[
\begin{aligned}
T_1&:\quad |I|<n,\\
T_2&:\quad |I|=n,\,\{1,m\}\subseteq I,\\
T_3&:\quad |I|=n,\quad\{1,m\}\nsubseteq I.
\end{aligned}
\]
Then
\[
\begin{aligned}
\cent{\Sn}&=\E(\Sn)\sqcup T_1\sqcup T_2\sqcup T_3,\\
\Mn(\Sn)&=\E(\Sn)\sqcup T_1\sqcup T_2.
\end{aligned}
\]
The inverse subsemigroups $\E(\Sn)$, $\E(\Sn)\cup T_1,$ and $\cent{\Sn}$ are all normal. Moreover, $T_2\ne\varnothing$ if and only if $m\ge n$, whereas $T_3\ne\varnothing$ if and only if $m>n$. When $m>n$, every element of $T_2$ has a conjugate in $T_3$.
\end{proposition}

\begin{proof}
Every singleton projection $E_{jj}$ belongs to $\Sn$. A matrix commuting with all of them is diagonal, and diagonal matrices commute with all idempotents in \eqref{eq:idempotents}. The partition and the margin formula follow. Diagonal products remain diagonal, their ranks cannot increase, and adjoints preserve rank; hence the three displayed sets are inverse subsemigroups. Conjugation preserves diagonal form and cannot increase rank. If a diagonal element of rank less than $n$ acquires trivial coefficient, its conjugate is idempotent; otherwise its rank remains less than $n$. This proves the stated normality assertions.

The nonemptiness conditions are immediate from the possible supports and $n\ge2$. Finally, the construction in \eqref{eq:conjugator} applies to every element of $T_2$ when $m>n$.
\end{proof}

\section{Marginal equivalence and inversion invariance}\label{sec:equivalence}

For an inverse word $\omega(x_1,\ldots,x_d)$, we recall that $a\sim_\omega b$ means that, for every $i$ and every choice of the variables,
\[
\begin{aligned}
\omega(x_1,\ldots,ax_i,\ldots,x_d)&=\omega(x_1,\ldots,bx_i,\ldots,x_d),\\
\omega(x_1,\ldots,x_i a,\ldots,x_d)&=\omega(x_1,\ldots,x_i b,\ldots,x_d).
\end{aligned}
\]
The marginal congruence $\equiv_\omega$ is the least inverse semigroup congruence containing $\sim_\omega$.

\begin{theorem}\label{thm:general-equivalence}
Let $S$ be an inverse semigroup, $\omega$ be an inverse semigroup word, and $\sim_\omega$ be the marginal equivalence relation associated with $\omega$. Then
\[
\sim_\omega=\equiv_\omega \quad\Longleftrightarrow\quad  M_\omega(S)\trianglelefteq S \ \text{and}\
\sim_\omega\text{is invariant under inversion}.
\]
\end{theorem}

\begin{proof}
If $\sim_\omega=\equiv_\omega$, then $\sim_\omega$ is a congruence and hence is invariant under inversion. For $a\in M_\omega(S)$ and $s\in S$, $a\sim_\omega a^0$ yields $sas^{-1}\sim_\omega sa^0s^{-1}$. Since $sa^0s^{-1}$ is idempotent, $sas^{-1}\in M_\omega(S)$. Thus $M_\omega(S)$ is self-conjugate and $M_\omega(S)\trianglelefteq S$.

Conversely, suppose that $M_\omega(S)\trianglelefteq S$ and $a\sim_\omega b$ implies $a^{-1}\sim_\omega b^{-1}$. Since $M_\omega(S)\subseteq\cent S$, Lawson in
\cite[Theorem~4.8]{LawsonPrimer} gives a congruence $\rho$ defined by $a\,\rho\,b$ if and only if $a^{-1}a=b^{-1}b$ and $ab^{-1}\in M_\omega(S)$. We show that $\rho=\sim_\omega$.

Let $a\,\rho\,b$ and $c=ab^{-1}\in M_\omega(S)$. Since $a^{-1}a=b^{-1}b$, we have
\[
cb=ab^{-1}b=a(a^{-1}a)=a.
\]
Since $c$ is completely regular, we also have
\[
c^0=c^{-1}c =ba^{-1}ab^{-1} =b(b^{-1}b)b^{-1} =bb^{-1}.
\]
It follows that $c^0b=bb^{-1}b=b$. Thus
\[
\begin{aligned}
\omega(\ldots,ax_i,\ldots) =\omega(\ldots,cbx_i,\ldots) =\omega(\ldots,c^0bx_i,\ldots) =\omega(\ldots,bx_i,\ldots).
\end{aligned}
\]
Moreover, for $t=x_icx_i^{-1}\in M_\omega(S)$, the inclusion $M_\omega(S)\subseteq\cent S$ gives $tx_i=x_ic$ and $t^0=x_ic^0x_i^{-1}$; hence $t^0x_ib=x_ib$. Thus
\[
\begin{aligned}
\omega(\ldots,x_ia,\ldots) =\omega(\ldots,tx_ib,\ldots) =\omega(\ldots,t^0x_ib,\ldots) = \omega(\ldots,x_ib,\ldots),
\end{aligned}
\]
and hence $a\sim_\omega b$.

Now let $a\sim_\omega b$. By inversion invariance, $a^{-1}\sim_\omega b^{-1}$. Substituting $b^{-1}x_i$ in the left identity for $a\sim_\omega b$ and $x_ia$ in the right identity for $a^{-1}\sim_\omega b^{-1}$ shows that $ab^{-1}\in M_\omega(S)$, with idempotent witnesses $bb^{-1}$ and $aa^{-1}$. Similarly, substituting $bx_i$ in the left identity for $a^{-1}\sim_\omega b^{-1}$ and $x_ia^{-1}$ in the right identity for $a\sim_\omega b$ shows that $a^{-1}b\in M_\omega(S)$, with witnesses $b^{-1}b$ and $a^{-1}a$. By \cite[Lemma~4]{AAK}, the two witnesses of each marginal element coincide. Thus $aa^{-1}=bb^{-1}$ and $a^{-1}a=b^{-1}b$, so $a\,\rho\,b$. Hence
$\sim_\omega=\rho$ is a congruence, and consequently $\sim_\omega=\equiv_\omega$.
\end{proof}

An inverse semigroup word $\omega$ is symmetric on $S$ if $\omega(x_1^{-1},\ldots,x_d^{-1})=\omega(x_1,\ldots,x_d)^{-1}$ for all $x_1,\ldots,x_d\in S$.

\begin{lemma}\label{lem:symmetric-inversion}
If $\omega$ is symmetric, then $\sim_\omega$ is invariant under inversion.
\end{lemma}
\begin{proof}
Suppose that $a\sim_\omega b$. Fix $i$ and arbitrary $x_1,\ldots,x_d\in S$. By symmetry and the right-insertion identity for $a\sim_\omega b$,
\[
\begin{aligned}
\omega(x_1,\ldots,a^{-1}x_i,\ldots,x_d) & =\omega(x_1^{-1},\ldots,x_i^{-1}a,\ldots,x_d^{-1})^{-1}\\
&=\omega(x_1^{-1},\ldots,x_i^{-1}b,\ldots,x_d^{-1})^{-1}=\omega(x_1,\ldots,b^{-1}x_i,\ldots,x_d).
\end{aligned}
\]
Similarly, symmetry and the left-insertion identity give
\[
\omega(x_1,\ldots,x_i a^{-1},\ldots,x_d)=\omega(x_1,\ldots,x_i b^{-1},\ldots,x_d).
\]
Thus $a^{-1}\sim_\omega b^{-1}$.
\end{proof}

\begin{corollary}\label{cor:symmetric-normality} 
If $\omega$ is symmetric, then 
\[
\sim_\omega=\equiv_\omega  \quad\Longleftrightarrow\quad  M_\omega(S)\trianglelefteq S.
\]
In particular, this equivalence holds for every power word $\omega(x)=x^n$.
\end{corollary}
\begin{proof}
This follows immediately from Theorem~\ref{thm:general-equivalence} and Lemma~\ref{lem:symmetric-inversion}. Indeed, $x^n$ is symmetric because $(x^{-1})^n=(x^n)^{-1}$.
\end{proof}

For a general inverse semigroup word, normality of the margin need not imply that marginal equivalence is a congruence. We show below that, for a suitable inverse semigroup word $\omega$, $S_{3,2}$ has a normal marginal inverse subsemigroup, but its marginal equivalence relation is not invariant under inversion.

\begin{proposition}\label{prop:omega-margin}
Let $S=S_{3,2}$ and $\omega(x,y)=x^2y^{-1}xy^{-1}x^{-1}$. Then
\[
M_\omega(S)=E(S)\cup\{-E_{11},-E_{22},-E_{33}\}.
\]
In particular, $M_\omega(S)$ is normal in $S$.
\end{proposition}

\begin{proof}
For $A\in M_\omega(S)$, complete regularity of marginal elements \cite[Lemma~4]{AAK} gives $AA^*=A^*A$. Since $AA^*=P_{R_A}$ and $A^*A=P_{D_A}$, we obtain $R_A=D_A$;
\eqref{eq:diagonal} then shows that $A$ is diagonal. The diagonal elements of $S$ are precisely those in $E(S)$, together with $-E_{ii}$ and $D_{ij}$ for $1\le i<j\le3$, where
\[
D_{ij}=-E_{ii}+E_{jj}
\qquad\text{and}\qquad
P_{ij}=D_{ij}D_{ij}^*=E_{ii}+E_{jj}.
\]
Since $E(S)\subseteq M_\omega(S)$, it remains to consider these non-idempotents.

For each $i$, we have $-E_{ii}\in M_\omega(S)$. Indeed, in either left or right insertion, replacing the idempotent witness $E_{ii}$ by $-E_{ii}$ contributes a factor $-1$ at every occurrence of the corresponding variable, including inverse occurrences. The word $\omega$ contains four occurrences of $x$ and two of $y$, counting inverse occurrences, so these signs cancel.

Let
\[
U=E_{21}+E_{33},\qquad
B=E_{11}+E_{32},\qquad
C=-E_{11}+E_{32}.
\]
Then $B^*=E_{11}+E_{23}$ and $C^*=-E_{11}+E_{23}$ satisfy
\[
\begin{aligned}
P_{12}B^*&=B^*=B^*P_{13},\\
D_{12}B^*&=C^*=B^*D_{13}.
\end{aligned}
\]
Since $\omega(U,B^*)=E_{32}=-\omega(U,C^*)$, the left-insertion test with $D_{12}$ and the right-insertion test with $D_{13}$ show that $D_{12},D_{13}\notin M_\omega(S)$. Moreover,
\[
\omega(B^*,P_{23}U)=E_{12}=-\omega(B^*,D_{23}U).
\]
Thus $D_{23}\notin M_\omega(S)$, proving the asserted formula
for $M_\omega(S)$.

Finally, for every $A\in S$ and every $i$, $A(-E_{ii})A^*$ is either $0$ or one of $-E_{11},-E_{22},-E_{33}$. Since $E(S)$ is full and self-conjugate, the displayed formula shows that $M_\omega(S)$ is normal.
\end{proof}

\begin{proposition}\label{prop:omega-inversion}
For $S$ and $\omega$ as in Proposition~\ref{prop:omega-margin}, the elements $B,C$ defined there satisfy
\[
B\sim_\omega C,\qquad B^*\not\sim_\omega C^*.
\]
Consequently, $\sim_\omega$ is not invariant under inversion.
\end{proposition}

\begin{proof}
Let $X,Y\in S$, with support maps $\phi_X,\phi_Y$. Write $\lambda_j,\mu_j\in\{\pm1\}$ for the unique nonzero entries in the $j$-th columns of $X$ and $Y$, for $j\in D_X$ and $j\in D_Y$, respectively. Every nonzero entry of $\omega(X,Y)$ corresponds to a chain
\[
a\xleftarrow{\phi_X}b\xleftarrow{\phi_X}c
 \xrightarrow{\phi_Y}d\xleftarrow{\phi_X}e
 \xrightarrow{\phi_Y}f\xrightarrow{\phi_X}g,
\]
so that
\[
\begin{gathered}
\phi_X(b)=a,\qquad \phi_X(c)=b,\qquad \phi_Y(c)=\phi_X(e)=d,\\
\phi_Y(e)=f,\qquad \phi_X(f)=g.
\end{gathered}
\]
The coefficient contributed by this chain is $\lambda_b\lambda_c\lambda_e\lambda_f\mu_c\mu_e$.

If $c=e$, then $b=d=f$, and this coefficient is $1$. If $c>e$, monotonicity of $\phi_X$ and $\phi_Y$ gives $b>d>f$. Since $b,c,e,f\in D_X$ and $\operatorname{rank}(X)\le2$, we have $D_X=\{e,c\}$, and hence $b=c$ and $f=e$. It follows that $c=3,e=1,d=2$, $\phi_X=\gamma$, and $\phi_Y=\beta$, where
\[
\beta(1)=1,\quad\beta(3)=2,\qquad
\gamma(1)=2,\quad\gamma(3)=3.
\]
Similarly, if $c<e$, then $b<d<f$, whence $b=c=1$, $e=f=3$, $d=2$, and $\phi_X=\beta,\ \phi_Y=\gamma$. In either case with $c\ne e$, we have $b=c$ and $e=f$. Thus the coefficient is independent of the coefficients of $X$, and can depend on those of $Y$ only if $\phi_Y\in\{\beta,\gamma\}$.

Let $\alpha$ be the increasing partial bijection with domain $\{1,2\}$, range $\{1,3\}$, and $\alpha(1)=1,\alpha(2)=3$. We have $\phi_B=\phi_C=\alpha$. For every $Z\in S$,
\[
R_{BZ}=R_{CZ}\subseteq\{1,3\},\qquad
D_{ZB}=D_{ZC}\subseteq\{1,2\}.
\]
Consequently, none of $\phi_{BZ},\phi_{CZ},\phi_{ZB},\phi_{ZC}$ belongs to $\{\beta,\gamma\}$. The coefficient calculation above now gives, for all $X,Y\in S$,
\[
\begin{aligned}
\omega(BX,Y)&=\omega(CX,Y),&
\omega(XB,Y)&=\omega(XC,Y),\\
\omega(X,BY)&=\omega(X,CY),&
\omega(X,YB)&=\omega(X,YC).
\end{aligned}
\]
Hence $B\sim_\omega C$.

On the other hand, $B^*=E_{11}+E_{23}$ and $C^*=-E_{11}+E_{23}$ satisfy
\[
\omega(U,B^*P_{13})=E_{32}=-\omega(U,C^*P_{13}),
\]
where $U=E_{21}+E_{33}$ is as in Proposition~\ref{prop:omega-margin}. Thus the right-insertion identity for the second variable fails, and $B^*\not\sim_\omega C^*$.
\end{proof}

Combining Propositions~\ref{prop:omega-margin} and \ref{prop:omega-inversion} with Theorem~\ref{thm:general-equivalence}, we obtain
\[
\sim_\omega\subsetneq\equiv_\omega
\quad\text{although}\quad
M_\omega(S)\trianglelefteq S
\]
for this example. Thus normality by itself does not imply that marginal equivalence is a congruence for an arbitrary inverse word.

\section{Two general observations on inverse words}\label{sec:generality}

We give two situations in which the bounds \eqref{eq:bounds} become equalities.

\begin{proposition}\label{prop:single}
Let $S$ be an inverse semigroup. If an inverse word $\omega$ contains some variable $x_i$ exactly once, either as $x_i$ or as $x_i^{-1}$, then $M_\omega(S)=\E(S).$ In particular, this margin is normal.
\end{proposition}

\begin{proof}
Let $a\in M_\omega(S)$, and set every variable equal to $a^0$. Replacing $x_i$ by $ax_i=a$ changes the value of the word from $a^0$ to $a$ or $a^{-1}$, according to the sign of its unique occurrence. Marginality forces $a=a^0$ or $a^{-1}=a^0$. Thus $a$ is idempotent. The reverse inclusion follows from \eqref{eq:bounds}.
\end{proof}

For the second observation, we recall the Munn representation \cite{Lawson}. Let $\mathord{\downarrow}e=\{f\in\E(S):f\le e\}$. For each $a\in S$, we define the partial bijection
\[
\theta_a\colon \mathord{\downarrow}(a^{-1}a) \longrightarrow \mathord{\downarrow}(aa^{-1}),
\qquad f\longmapsto afa^{-1}.
\]
The assignment $\theta\colon a\mapsto\theta_a$ defines an inverse semigroup homomorphism from $S$ to the symmetric inverse monoid on $\E(S)$. This homomorphism separates idempotents.

\begin{lemma}\label{lem:centralizer}
An element $a\in S$ belongs to $\cent S$ if and only if $\theta_a$ is a partial identity. For such $a$,
\[
aa^{-1}=a^{-1}a=a^0,\qquad \theta_a=\theta_{a^0}.
\]
Consequently $\cent S$ is a normal inverse subsemigroup.
\end{lemma}

\begin{proof}
Suppose that $a$ centralizes all idempotents, and let $d=a^{-1}a$, $r=aa^{-1}$. The equalities $ar=ra=a$ and $da=ad=a$ give $dr=d$ and $dr=r$, respectively. Thus $d=r=a^0$. For $f\le a^0$,
\[
\theta_a(f)=afa^{-1}=faa^{-1}=f,
\]
proving $\theta_a=\theta_{a^0}$.

Conversely, if $\theta_a$ is a partial identity, its domain and range coincide, so $a^{-1}a=aa^{-1}=e$. For arbitrary $f\in\E(S)$, apply $\theta_a$ to $ef$ to get $aef a^{-1}=ef$. Multiplication on the right by $a$ gives $af=efa=fa$. This proves the equivalence. The inverse image under $\theta$ of the idempotents is full, closed under multiplication and inversion, and self-conjugate. The last assertion follows.
\end{proof}

\begin{proposition}\label{prop:idempotent-values}
Let $S$ be an inverse semigroup and let $\omega$ be an inverse word such that every value of $\omega$ in $S$ is idempotent. Then
$$
M_\omega(S)=\cent S.
$$
In particular, this margin is normal.
\end{proposition}

\begin{proof}
By \eqref{eq:bounds}, it is enough to prove the reverse inclusion. Take $a\in\cent S$. Lemma~\ref{lem:centralizer} gives $\theta_a=\theta_{a^0}$, and therefore
\[
\theta_{ax_i}=\theta_{a^0x_i}
\quad\text{and}\quad
\theta_{x_ia}=\theta_{x_ia^0}.
\]
Since $\theta$ respects products and inverses, the two word values
\[
\begin{aligned}
p&=\omega(x_1,\ldots,ax_i,\ldots,x_d),\\
q&=\omega(x_1,\ldots,a^0x_i,\ldots,x_d)
\end{aligned}
\]
satisfy $\theta_p=\theta_q$. Both are idempotent by assumption, so the idempotent-separating property implies $p=q$. This proves \eqref{eq:def-left} with $e_i=a^0$. The same argument with right insertion proves \eqref{eq:def-right}, and hence $a\in M_\omega(S)$.
\end{proof}

\section*{Acknowledgements}
The authors thank zbMATH for its review of \cite{AAK}, which drew their attention to Problem \ref{AAK-Problem11}. They also acknowledge the use of AI tools to assist in revising the manuscript. This work was partially supported by the National Natural Science Foundation of China (Grant No. 12001254) and by Lanzhou University (No. 561120205).

\end{document}